\documentclass{autart}
\usepackage{mathtools}
\usepackage{amstext}
\usepackage{amssymb}

\usepackage{bm}
\usepackage{color}

\renewcommand{\S}{\mathcal{S}}
\newcommand{\I}{\mathcal{I}}
\newcommand{\G}{\mathcal{G}}
\newcommand{\V}{\mathcal{V}}
\newcommand{\E}{\mathcal{E}}
\newcommand{\x}{\mathbf{x}}
\renewcommand{\u}{\mathbf{u}}
\newcommand{\w}{\mathbf{w}}
\renewcommand{\c}{\mathbf{c}}
\newcommand{\boldb}{\bm{\beta}}
\newcommand{\boldg}{\bm{\gamma}}
\newcommand{\tL}{\mathcal{L}}
\newcommand{\bin}[2]{b^{#1}_{#2}}
\newcommand{\eigsum}{\sigma}
\usepackage{flushend}
\newtheorem{rmk}[thm]{Remark}
\allowdisplaybreaks

\usepackage{graphicx}      
\usepackage[labelformat=simple]{subcaption}

\usepackage{natbib}
\begin{document}
\begin{frontmatter}

\title{A Dynamical Approach to Efficient Eigenvalue Estimation in General Multiagent Networks} 

\thanks[footnoteinfo]{This work was supported, in part, by the National Science Foundation, grants CAREER-ECCS-1651433 and III-200884556.}

\author[First]{Mikhail Hayhoe}
\author[First]{Francisco Barreras}
\author[First]{Victor M. Preciado} 

\address[First]{University of Pennsylvania, 
   Philadelphia, PA 19143 USA (e-mails: mhayhoe@seas.upenn.edu, fbarrer@sas.upenn.edu, preciado@seas.upenn.edu).}

\begin{abstract}                
We propose a method to efficiently estimate the eigenvalues of any arbitrary (potentially weighted and/or directed) network of interacting dynamical agents from dynamical observations. These observations are discrete, temporal measurements about the evolution of the outputs of a subset of agents (potentially one) during a finite time horizon; notably, we do not require knowledge of which agents are contributing to our measurements. We propose an efficient algorithm to exactly recover the (potentially complex) eigenvalues corresponding to network modes that are observable from the output measurements. The length of the sequence of measurements required by our method to generate a full reconstruction of the observable eigenvalue spectrum is, at most, twice the number of agents in the network, but smaller in practice. The proposed technique can be applied to networks of multiagent systems with arbitrary dynamics in both continuous- and discrete-time. Finally, we illustrate our results with numerical simulations.
\end{abstract}

\begin{keyword}
Multiagent networks; eigenvalue estimation; sparse estimation; spectral identification; Laplacian matrix
\end{keyword}

\end{frontmatter}

\section{Introduction}

The spectra of matrices describing the structure of a network of interacting dynamical agents provide a wealth of global information about the network structure and function; see, e.g., \citet{Fiedler73, mohar91, merris94, Chung97, Preciado08, mesbahi10, Bullo19}, and references therein.
A particular example of interest is the Laplacian spectrum, which finds applications in multiagent coordination problems \citep{jadbabaie03, olfati07}, synchronization of oscillators \citep{pecora98, dorfler13}, neuroscience \citep{Becker18}, biology \citep{Palsson06}, as well as several graph-theoretical problems, such as finding cuts \citep[see][]{shi00} or communities \citep[see][]{Luxburg07} in graphs, among many others, as illustrated in \cite{mohar97}. Beyond the Laplacian eigenvalues, the spectrum of the adjacency matrix of a network is relevant in the analysis of, for example, epidemic processes \citep{Nowzari16}. Furthermore, the eigenvalues of the normalized Laplacian are relevant in the analysis of diffusion processes, random walks over graphs, or discrete-time consensus dynamics \citep{Chung97}.

Due to its practical importance, numerous methods have been proposed to estimate the eigenvalues of a network of dynamical agents.
For example, \cite{kempe08} proposed a distributed algorithm based on orthogonal iteration \citep[see][]{golub13} for computing higher-dimensional invariant subspaces.
In the control literature, \citet{Franceschelli13} define local interaction rules between agents such that the network response is a superposition of sinusoids oscillating at frequencies related to the Laplacian eigenvalues; however, this approach imposes a particular dynamics on the agents in the network, which is unrealistic in many scenarios.
\citet{Aragues14} proposed a distributed algorithm based on the power iteration for computing upper and lower bounds on the algebraic connectivity (i.e., the second smallest Laplacian eigenvalue).
\citet{Leonardos20} proposed a distributed continuous-time dynamics over manifolds to compute the largest (or smallest) eigenvalues and eigenvectors of any graph.
An approach by \citet{Kibangou15} uses consensus optimization to deduce the spectrum of the Laplacian, but this requires a consensus algorithm to be run on the network separately from the dynamics.
Using the Koopman operator, it has been shown that the spectrum of the Laplacian may be recovered using sparse local measurements, see \cite{Mauroy17,Mesbahi19}; 
unfortunately, these methods require the system to be reset to known initial conditions multiple times or for full observability of agents' states.

We find in the literature several works more closely related to the techniques used in this paper. For example, a classical approach known as \emph{Prony's method} can be used to estimate the parameters of a uniformly sampled superposition of complex exponentials, which can be used for spectral estimation and deconvolution, among other problems \citep[see][]{Potts10,Kunis16}. In contrast to our approach, Prony's method only applies to symmetric matrices; hence, it can only be applied for the spectra lreconstruction of undirected networks. Also related to our work we find the \emph{Newton-Girard equations} \citep[see, e.g.,][]{Herstein06} which allow us to recover eigenvalues by analyzing symmetric polynomials of the traces of powers of the matrix. However, computing the traces of powers of matrices is computationally expensive and requires a large amount of (centralized) data, which may not be feasible to collect in many applications. Using local structural information, \cite{Preciado13-2} computed the traces of powers of a graph matrix to derive bounds on spectral properties of practical importance, such as the spectral radius. A related method uses the classical moment problem from probability theory to analyze the spectrum of a graph by counting walks in graphs, as in \cite{Preciado13, Chen20, Barreras19}.

In this paper we present an approach to efficiently estimate the eigenvalues of any graph matrix, such as the Laplacian, corresponding to an unknown network of multiagent systems using only a single temporal sequence of (potentially sparse) output measurements. The network structure may be weighted and/or directed, and may include multi-edges and self-loops. The temporal sequence measurements used in our spectral estimation algorithm can correspond to the output signal of a single agent, or to any weighted linear combination of outputs from a collection of agents; notably, our method requires no knowledge of which agents contribute to the measurements, nor does it require prior knowledge of the network topology or its initial condition. Moreover, the length of the sequence of measurements required is, at most, twice the number of agents in the network, but fewer in practice. Our approach allows for the estimation of all complex eigenvalues associated with observable network modes, regardless of the (unknown) network structure. The proposed approach requires no tuning of parameters, and may be applied in both discrete- and continuous-time to general multi-agent systems.

The remainder of this paper is structured as follows. We outline background and notation in Section~\ref{sec:bg}. We introduce our approach on the particular case of discrete-time Laplacian dynamics in Section~\ref{sec:dt_laplacian}. In Section~\ref{sec:dt} we present our results for discrete-time systems and in Section~\ref{sec:ct}, we describe our results in the continuous-time case. Section~\ref{sec:sims} illustrates our results via simulations in a variety of systems, and Section~\ref{sec:conc} concludes the paper.

\section{Background and Notation}\label{sec:bg}

\begin{table}[ht]
    \centering
    \begin{tabular}{c|l}
        Symbol & Meaning \\
        \hline
        $I_n$ & $n\times n$ identity matrix \\
        $\mathbb{R}$ & set of real numbers \\
        $\mathbb{N}$ & set of natural numbers \\
        $\mathbf{e}_i$ & $i$-th vector in the canonical basis of $\mathbb{R}^n$ \\
        $\V$ & node set, $\V = \{1,\ldots,n\}$ \\
        $\E$ & edge set, $\E \subseteq \V\times\V$ \\
        $\G = (\V,\E) $ & graph with node set $\V$ and edge set $\E$\\
        $\otimes$ & Kronecker product \\
        $\oplus$ & Direct sum \\
        $\!\sigma(X) \!\coloneqq\! \{\lambda_i\}_{i=1}^n\!$ & eigenvalue spectrum of matrix $X$ \\
        $A(\G)$ & adjacency matrix of $\G$, $[A]_{ij} \neq 0 \Rightarrow (i,j)\in\E$ \\
        $D(\G)$ & degree matrix of $\G$, $[D]_{ii} = \sum_{j=1}^n [A]_{ij}$ \\
    \end{tabular}
    \label{tab:notation}
\end{table}

Throughout this paper we use lower-case letters for scalars, lower-case bold letters for vectors, upper-case letters for matrices, and calligraphic letters for sets.

A \emph{directed} graph $\G = (\V,\E)$ has node set $\V$ and edge set $\E$, where $(i,j)\in\E$ means node $i$ has an edge pointed toward node $j$. The graph $\G$ may have self-loops, may have (possibly negative) edge weights, and may contain multi-edges. 


\section{Discrete-Time Laplacian Dynamics}\label{sec:dt_laplacian}

We begin our exposition with a simple, undirected network of single integrators following a discrete-time (DT) Laplacian dynamics. In this context, we will introduce a methodology to estimate the eigenvalues of the Laplacian matrix from a finite sequence of output measurements; for full details of this case, see~\cite{Hayhoe19}. In Section~\ref{sec:dt}, we will extend this result to more general directed networks of discrete-time agents, and will study the continuous-time (CT) case in Section~\ref{sec:ct}.

Consider the following discrete-time dynamics:
\begin{align}\begin{split}\label{eq:lap_integrator}
\mathbf{x}\left[k+1\right] & = \tL\mathbf{x}\left[k\right],\;\mathbf{x}\left[0\right]=\mathbf{x}_{0},\\
y\left[k\right] & =\mathbf{c}^{\intercal}\mathbf{x}\left[k\right],
\end{split}
\end{align}
where $\tL \coloneqq D(\G)^{-1}A(\G)$ is the normalized Laplacian matrix of an unknown undirected graph $\G$, $k\in\mathcal{\mathbb{N}}$, and $\mathbf{c},\mathbf{x}_{0}$
are arbitrary (possibly unknown) vectors in $\mathbb{R}^{n}$. For example, we may have
$\mathbf{c}=\mathbf{e}_{i}$ when we only observe the state of agent
$i$, or $\mathbf{c}=\sum_{i\in\mathcal{S}}\beta_i\mathbf{e}_{i}$
when we observe the weighted sum of the states of a subset $\mathcal{S} \subseteq \mathcal{V}$ of
agents. In what follows we propose an efficient algorithm to recover the eigenvalues of the normalized Laplacian matrix $\tL$ from the output sequence $y[0],y[1],\ldots,y[2n-1]$.

The normalized Laplacian $\tL$ of an undirected graph is always diagonalizable with real eigenvalues $\lambda_1,\ldots,\lambda_n \in \mathbb{R}$ \citep[see][]{Chung97}. Denoting by $\mathbf{u}_i$ and $\mathbf{w}_i$ the (unknown) right and left eigenvectors corresponding to the eigenvalue $\lambda_i$, we have that $ \tL = U\Lambda W$, where $\Lambda \coloneqq \text{diag}(\lambda_1,\ldots,\lambda_n)$, $U \coloneqq [\u_1,\ldots,\u_n]$, and  $W \coloneqq [\w_1^\intercal;\cdots;\w_n^\intercal] = U^{-1}$; hence,
\begin{align}\label{eq:Laplacian_output}
y\left[k\right] & =\mathbf{c}^{\intercal}\tL^{k}\mathbf{x}_{0} =\left(\mathbf{c}^{\intercal}U\right)\Lambda^{k}\left(W\mathbf{x}_{0}\right) =\sum_{i=1}^{n}\omega_{i}\lambda_{i}^{k},
\end{align}
where the weights $\omega_i$ are given by
\begin{equation}
\omega_{i} \coloneqq \left[\mathbf{c}^{\intercal}U\right]_{i}\left[W\mathbf{x}_{0}\right]_{i}=\mathbf{c}^{\intercal}\mathbf{u}_{i}\mathbf{w}_{i}^{\intercal}\mathbf{x}_{0}.\label{eq:weights}
\end{equation}
Notice that it is possible for $\omega_{i}=0$ whenever $\mathbf{c}^{\intercal}\mathbf{u}_{i}=0$ or $\mathbf{w}_{i}^{\intercal}\mathbf{x}_{0}=0$. If $\omega_i = 0$ for some index $i$, then the $i$-th eigenvalue $\lambda_i$ does not influence the output $y[k]$ in~\eqref{eq:Laplacian_output}; consequently, we will not be able to estimate $\lambda_i$ from a sequence of outputs. However, if $\x_0$ is randomly generated, then almost surely $\w_i^\intercal\x_0 \neq 0$; hence, it is possible that $\omega_i = 0$ only for those eigenvalues $\lambda_i$ for which $\mathbf{c}^{\intercal}\mathbf{u}_{i}=0$. According to the Popov-Belevitch-Hautus (PBH) test \citep[see][]{hespanha18}, those eigenvalues corresponding to unobservable eigenmodes of the Laplacian dynamics will be those for which $\omega_i=0$ and it will be impossible to recover them from our observations.
Furthermore, we can have repeated eigenvalues that would not impact
the output $y\left[k\right]$, whenever $\sum_{j\colon\lambda_{j}=\lambda_{i}}\omega_{j}=0$.
Defining the constant $\mathbf{w}_{j}^{\intercal}\mathbf{x}_{0}=\alpha_{j}$
(which will be different than zero almost surely), this condition
is equivalent to 
\[
\sum_{j\colon\lambda_{j}=\lambda_{i}}\omega_{j}=\mathbf{c}^{\intercal}\sum_{j\colon\lambda_{j}=\lambda_{i}}\alpha_{j}\mathbf{u}_{j}=\mathbf{c}^{\intercal}\mathbf{u}^{\left(\lambda_{i}\right)}=0
\]
where $\mathbf{u}^{\left(\lambda_{i}\right)}=\sum_{j\colon\lambda_{j}=\lambda_{i}}\alpha_{j}\mathbf{u}_{j}$.
Note that $\mathbf{u}^{\left(\lambda_{i}\right)}$ is in the eigenspace
of the eigenvalue $\lambda_{i}$; hence, according to the PBH test,
$\mathbf{c}^{\intercal}\mathbf{u}^{\left(\lambda_{i}\right)}=0$ implies
that $\lambda_{i}$ corresponds to an unobservable eigenmode (almost surely). We will denote by $\S_{\tL}$ the set of eigenvalues of $\tL$ corresponding to observable eigenmodes of the pair $(\tL,\mathbf{c}^\intercal)$.

In the theorem below, we describe a methodology to efficiently reconstruct the observable eigenvalues $\lambda_i \in \S_{\tL}$ from a finite sequence of output observations.

\begin{thm}\label{thm:id}
Given the sequence of observations $\left(y\left[k\right]\right)_{k=0}^{2n-1}$ from the system in \eqref{eq:dt_integrator}, define the following Hankel matrix
\begin{align}\label{eq:hankel_Y}
Y \coloneqq \left[\begin{array}{cccc}
y[0] & y[1] & \cdots & y[n-1]\\
y[1] & y[2] & \cdots & y[n]\\
\vdots & \vdots & \ddots & \vdots\\
y[n-1] & y[n] & \cdots & y[2n-2]
\end{array}\right].
\end{align}
The rank of the Hankel matrix $Y$ satisfies
\[
	r \coloneqq \text{rk}(Y) = |\S_{\tL}| \leq n.
\]
The observable eigenvalues of $\tL$ are roots of the polynomial
\[
p_{\tL}\left(x\right)=x^{r}+\alpha_{r-1}x^{r-1}+\cdots+\alpha_{1}x+\alpha_{0},
\]
where the coefficients $\alpha_{0},\ldots,\alpha_{r-1}$ are given
by
\[
\!\left[\!\begin{array}{c}
\alpha_{0}\\
\alpha_{1}\\
\vdots\\
\alpha_{r-1}
\end{array}\!\right]\!=\!-\!\left[\!\begin{array}{cccc}
y[0] & y[1] & \cdots & y[r-1]\\
y[1] & y[2] & \cdots & y[r]\\
\vdots & \vdots & \ddots & \vdots\\
y[r-1] & y[r] & \cdots & y[2r-2]
\end{array}\right]^{-1}\!\left[\!\begin{array}{c}
y[r]\\
y[r+1]\\
\vdots\\
y[2r-1]
\!\end{array}\right]\!.
\]
\end{thm}
\begin{pf}
See~\cite{Hayhoe19}.
\end{pf}

In what follows, we will extend this result to any arbitrary (possibly weighted and/or directed) network, in both discrete- and continuous-time.

\section{Spectral Estimation for Discrete-Time Dynamics}\label{sec:dt}

Let $G$ be any graph matrix whose sparsity pattern describes the connections of an arbitrary (unknown) graph $\G$ with $n$ nodes. The graph $\G$ may be directed, may have self-loops, and may be weighted. Consider the discrete-time dynamics of a collection of single integrators,
\begin{align}\begin{split}\label{eq:dt_integrator}
\mathbf{x}\left[k+1\right] & = G\mathbf{x}\left[k\right], \; \mathbf{x}\left[0\right]=\mathbf{x}_{0},\\
y\left[k\right] & =\mathbf{c}^{\intercal}\mathbf{x}\left[k\right],
\end{split}
\end{align}
where $k\in\mathcal{\mathbb{N}}$, and $\mathbf{c},\mathbf{x}_{0}$
are arbitrary (possibly unknown) vectors in $\mathbb{R}^{n}$. We may view our approach as a decentralized estimation problem when $\mathbf{c}=\mathbf{e}_{i}$, wherein agent $i$ is attempting to estimate the eigenvalues of $\G$ by observing its own output. More generally, we may observe the weighted sum of the states of a subset $\mathcal{S} \subseteq \mathcal{V}$ of agents; hence, $\mathbf{c}=\sum_{i\in\mathcal{S}}\beta_i\mathbf{e}_{i}$, which corresponds to a group of agents collectively estimating the spectrum of $\G$ using a weighted linear combination of their outputs using (possibly unknown) weights $\{\beta_i\}_{i\in\S}$.

To extend the result in Section~\ref{sec:dt_laplacian} to more general (possibly weighted and/or directed) dynamics, we start by defining the Jordan decomposition of $G$ as
\[
G = VJV^{-1}=V\left[\begin{array}{cccc}
J_{1} & 0 & \cdots & 0\\
0 & J_{2} & \cdots & 0\\
\vdots & \vdots & \ddots & \vdots\\
0 & 0 & \cdots & J_{d}
\end{array}\right]V^{-1},
\]
where $J_{i},~i\in\left\{ 1,\ldots,d\right\} $, is the $m_i \times m_i$ Jordan block associated with the $i$-th eigenvalue $\lambda_{i}$. Note that there may be
multiple Jordan blocks associated with a single eigenvalue; hence,
it may be that $\lambda_{i}=\lambda_{j}$ for some $i,j\in\left\{ 1,\ldots,d\right\} $. We thus also define the largest block size for each distinct eigenvalue $\lambda_i$ as $\hat{m}_i \coloneqq \max_{j : \lambda_j = \lambda_i} m_j$. Taking powers of the matrix $G$, we obtain $G^{k}=\left(VJV^{-1}\right)^{k}=VJ^{k}V^{-1}$,
where the $m_i \times m_i$ Jordan block raised to the power $k$, $J_{i}^{k}$, is the upper-triangular matrix
\begin{align}\label{eq:Jordan_exp}
J_{i}^{k} \!=\! \left[\!\begin{array}{ccccc}
\!\lambda_{i}^{k} & {k \choose 1}\lambda_{i}^{k-1} & {k \choose 2}\lambda_{i}^{k-2} & \!\cdots\! & {k \choose m_{i}-1}\lambda_{i}^{k-(m_{i}-1)}\! \\
 & \lambda_{i}^{k} & {k \choose 1}\lambda_{i}^{k-1} & \!\cdots\! & {k \choose m_{i}-2}\lambda_{i}^{k-(m_{i}-2)}\! \\
 &  & \ddots &  & \vdots\\
 &  &  & \lambda_i^k & {k \choose 1}\lambda_i^{k-1} \\
 &  &  &  & \lambda_{i}^{k}
\end{array}\!\right]\!.
\end{align}
For $i\in\left\{ 1,\ldots,d\right\} $, let
$\left[\mathbf{c}^{\intercal}V\right]_{i}$ and $\left[V^{-1}\mathbf{x}_{0}\right]_{i}$
denote the $m_i$-dimensional $i$-th blocks of $\mathbf{c}^{\intercal}V$ and $V^{-1}\mathbf{x}_{0}$,
respectively, associated with Jordan block matrix $J_{i}$. Hence, for any graph matrix $G$, the observations from our system~\eqref{eq:dt_integrator} can be written as
\begin{align}
y\left[k\right] & =\left(\mathbf{c}^{\intercal}V\right)J^{k}\left(V^{-1}\mathbf{x}_{0}\right)\label{eq:obs}\\
 & =\sum_{i=1}^{d}\left[\mathbf{c}^{\intercal}V\right]_{i}J_{i}^{k}\left[V^{-1}\mathbf{x}_{0}\right]_{i}\nonumber \\
 & =\sum_{i=1}^{d}\sum_{s=0}^{m_{i}-1}\omega_{i}^{(s)}{k \choose s}\lambda_{i}^{k-s},\nonumber 
\end{align}
where, for $s \!\in\! \{0,\ldots,m_i-1\}$, the weights $\omega_{i}^{(s)}$ are defined as
\begin{equation}
\omega_{i}^{(s)}\coloneqq\sum_{l=s+1}^{m_{i}}\left[\mathbf{c}^{\intercal}V\right]_{i,l-s}\left[V^{-1}\mathbf{x}_{0}\right]_{i,l},\label{eq:weights}
\end{equation}
with $\left[\mathbf{c}^{\intercal}V\right]_{i,l}$ and $\left[V^{-1}\mathbf{x}_{0}\right]_{i,l}$,
$l\in\left\{ 1,\ldots,m_{i}\right\} $ being the $l$-th components of
$\left[\mathbf{c}^{\intercal}V\right]_{i}$ and $\left[V^{-1}\mathbf{x}_{0}\right]_{i}$,
respectively. Finally, define the total weights corresponding to each unique eigenvalue as
\begin{align}\label{eq:total_weights}
	\bar{\omega}_{i}^{(s)} \coloneqq \sum_{j:\lambda_{j}=\lambda_{i}}\omega_{j}^{(s)}.
\end{align}
In general it is possible that $\bar{\omega}_{i}^{(s)}=0$, which will make it impossible to recover the eigenvalue $\lambda_i$. According to the PBH test, these eigenvalues correspond to unobservable eigenmodes of the pair $(G, \mathbf{c}^\intercal)$. We denote the set of observable eigenvalues by
\begin{align}\label{eq:S}
	\S_{G} \coloneqq \left\{ \lambda_{i}\in\sigma\left(G\right) \colon \exists s \text{ s.t. } \bar{\omega}_i^{(s)} \neq 0\right\}.
\end{align}
For an eigenvalue $\lambda_i \in \S_G$, we define
\begin{align} \label{eq:m_tilde}
	\tilde{m}_i \coloneqq 1 + \max\left\{s = 0,\ldots,\hat{m}_i-1 : \bar{\omega}_i^{(s)} \neq 0\right\},
\end{align}
and denote the set of indices corresponding to unique\footnote{Unique eigenvalues refers to unique values, i.e., the eigenvalues ignoring multiplicity.} observable eigenvalues as $\I \coloneqq \left\{ i\in\{1,\ldots,n\} : \lambda_i \in \S_G \right\}$. We can thus rewrite the observations from~\eqref{eq:obs} as
\begin{align}\label{eq:obs_total}
	y[k] = \sum_{i\in\I} \sum_{s=0}^{\hat{m}_i-1} \bar{\omega}_i^{(s)}\binom{k}{s}\lambda_i^{k-s}.
\end{align}
%

In what follows, we will propose a computationally efficient methodology to recover the eigenvalues in $\S_G$ using the output sequence $\left(y\left[k\right]\right)_{k=0}^{2n-1}$. Towards that goal, we define the Hankel matrix of observations
\begin{equation}
H \coloneqq \left[\begin{array}{cccc}
y[0] & y[1] & \cdots & y[n-1]\\
y[1] & y[2] & \cdots & y[n]\\
\vdots & \vdots & \ddots & \vdots\\
y[n-1] & y[n] & \cdots & y[2n-2]
\end{array}\right].\label{eq:hankel}
\end{equation}
The following result relates the rank of this matrix to the largest observable Jordan blocks of $G$.

\begin{lem}
\label{lem:rank}The rank of $H$ in \eqref{eq:hankel} satisfies
\[
\text{rk}(H) = \sum_{i \in \I} \tilde{m}_i,
\]
where $\tilde{m}_i$ is defined in~\eqref{eq:m_tilde}.
\end{lem}
\begin{pf}
See Appendix~\ref{app:lem_rank_proof}.
\end{pf}

With this Lemma in hand, we present our main result on estimating the observable eigenvalues of the pair $(G, \mathbf{c}^\intercal)$.

\begin{thm}\label{thm:id}
Given the sequence of observations $\left(y\left[k\right]\right)_{k=0}^{2n-1}$ from the discrete-time system in \eqref{eq:dt_integrator}, consider the matrix $H$ from~\eqref{eq:hankel} and denote its rank by $r$. The observable eigenvalues are roots of the polynomial
\[
p_{G}\left(x\right)=x^{r}+\alpha_{r-1}x^{r-1}+\cdots+\alpha_{1}x+\alpha_{0},
\]
where the coefficients $\alpha_{0},\ldots,\alpha_{r-1}$ are given
by
\[
\left[\!\begin{array}{c}
\alpha_{0}\\
\alpha_{1}\\
\vdots\\
\alpha_{r-1}
\end{array}\!\right]\!=\!-\!\left[\begin{array}{cccc}
y[0] & y[1] & \cdots & y[r-1]\\
y[1] & y[2] & \cdots & y[r]\\
\vdots & \vdots & \ddots & \vdots\\
\!y[r-1] & y[r] & \cdots & y[2r-2]\!
\end{array}\right]^{\!-1}\!\left[\begin{array}{c}
y[r]\\
y[r+1]\\
\vdots\\
\!y[2r-1]\!
\end{array}\right]\!.
\]
Moreover, $\lambda_i \!\in\! \S_G$ is a root of $p_G(x)$ with multiplicity $\tilde{m}_i$.
\end{thm}
\begin{pf}
See Appendix~\ref{app:thm_proof}.
\end{pf}

\begin{rmk}
While Theorem~\ref{thm:id} makes use of $2n$ observations $\left(y\left[k\right]\right)_{k=0}^{2n-1}$, in practice, fewer observations may be required. To illustrate this, consider an online setting where the output measurements are taken sequentially, one at a time. We can build a $k\times k$ Hankel matrix using the first $2k-1$ observations from the system, $y[0],y[1],\ldots,y[2k-2]$, and check its rank. If it is full rank we continue taking measurements. By the structure of the Hankel matrix, if the rank does not grow after including measurements $y[2k-1]$ and $y[2k]$ then it must be that there exists some non-trivial $\alpha_0,\ldots,\alpha_{k-1}$ such that $y[k] = \alpha_0 y[0] + \cdots \alpha_{k-1}y[k-1]$. By definition of the observations in~\eqref{eq:obs}, we thus have
\begin{align*}
\c^\intercal VJ^{k}V^{-1}\x_0 &= \sum_{s=0}^{k-1}\alpha_s\c^\intercal VJ^sV^{-1}\x_0 \\
	&= \c^\intercal V \left(\sum_{s=0}^{k-1}\alpha_sJ^s \right)V^{-1}\x_0.
\end{align*}
If all eigenmodes are observable and the initial condition $\x_0$ is random, then by Cayley-Hamilton theorem we must have $k=n$ (almost surely). If $l$ of the eigenmodes are unobservable (including multiplicities), then necessarily $k \geq n - l$, since at most $l$ entries could be zeroed by the vector $\c^\intercal V$. Thus, the rank of the Hankel matrix will stop growing once $k = r$ (almost surely), i.e., after we have collected enough measurements to recover all distinct eigenvalues corresponding to observable eigenmodes.
\end{rmk}

\subsection{Network of Identical Discrete-Time Agents}\label{subsec:dt_multi}

In many applications, the network of interest will consist of agents with more general dynamics beyond single integrators. With this in mind, we consider a network of $n$ agents where each agent follows the dynamics
$\mathbf{x}_{i}\left[k+1\right]=A\mathbf{x}_i\left[k\right]+\mathbf{u}_i\left[k\right]$, where\textbf{ $\mathbf{x}_{i}$} is a $d$-dimensional vector of states, $A$ is a known $d\times d$ state transition matrix, and $\mathbf{u}_i\left[k\right]$ is an input consisting of a linear combination of the states of
the neighboring agents of $i$. Assuming that all agents start with an arbitrary initial condition $\boldb$ weighted by $x_{0i}$, and the output of agent $i$ is $\boldg^\intercal\x_i[k]$  weighted by $c_i$, we obtain the following network dynamics:
\begin{align}
\begin{split}\label{eq:dt_multi}
\mathbf{x}_{i}\left[k+1\right] &=A\mathbf{x}_i\left[k\right]+\sum_{j=1}^{n}g_{ij}\mathbf{x}_{j}\left[k\right],\;\mathbf{x}_{i}\left[0\right]=x_{0i}\boldb,\\
y\left[k\right] & =\sum_{i=1}^{n}c_{i}\boldg^{\intercal}\mathbf{x}_{i}\left[k\right],
\end{split}
\end{align}
where $g_{ij} = [G]_{ij}$, $c_i = [\mathbf{c}]_i$, $x_{0i} = [\x_0]_i$. Stacking the vectors of states in a large vector $\mathbf{x}=\left(\mathbf{x}_{1}^{\intercal},\ldots,\mathbf{x}_{n}^{\intercal}\right)^{\intercal}$,
the dynamics can be written as
\begin{align*}
\mathbf{x}\left[k+1\right] & =\left(I_{n}\otimes A+ G\otimes I_{d}\right)\mathbf{x}\left[k\right],\;\mathbf{x}\left[0\right]=\mathbf{x}_0\otimes\boldb,\\
y\left[k\right] & =\left(\mathbf{c}\otimes\boldg\right)^{\intercal}\mathbf{x}\left[k\right].
\end{align*}
We assume the state matrix $A$ as well as the vectors of individual initial condition $\boldb$ and observation $\boldg$ are known, but the graph matrix $G$ and weighting vectors for initial conditions $\x_0$ and observations $\mathbf{c}$ are unknown. Our aim is to estimate the observable eigenvalues of $G$ from a finite sequence of outputs. This result is stated in the following theorem.

\begin{thm}\label{thm:id_multi}
Given the sequence of observations $\left(y\left[k\right]\right)_{k=0}^{2n-1}$ from the system in \eqref{eq:dt_multi}, consider the Hankel matrix $H$ defined in \eqref{eq:hankel} and denote its rank by $r$. The weighted sums of eigenvalues $\eigsum_k \coloneqq \sum_{i=1}^d\sum_{s=0}^{m_{i}-1}\omega_{i}^{(s)}{k \choose s}\lambda_{i}^{k-s}$  satisfy the following equality:
\begin{align*}
\!\begin{bmatrix}
\eigsum_{0}\\
\eigsum_{1}\\
\vdots\\
\!\eigsum_{2r-1}
\end{bmatrix} \!=\! \begin{bmatrix}
\bin{0}{0}\nu_{0} & 0 & \cdots & 0\\
\bin{1}{0}\nu_{1} & \bin{1}{1}\nu_{0} & \cdots & 0\\
\vdots & \vdots & \ddots & \vdots\\
\!\bin{2r-1}{0}\nu_{2r-1} & \bin{2r-1}{1}\nu_{2r-2} \!&\! \cdots \!&\! \bin{2r-1}{2r-1}\nu_{0}\!
\end{bmatrix}^{\!-1\!}
\begin{bmatrix}
y_{0}\\
y_{1}\\
\vdots\\
\!y_{2r-1}
\end{bmatrix}
\end{align*}
where $\nu_{k-s} \coloneqq \boldg^{\intercal}A^{k-s}\boldb$, $\bin{k}{s} \coloneqq \binom{k}{s}$, and the matrix is invertible when $\boldg^\intercal\boldb \neq 0$. Then, the observable eigenvalues of $G$ are roots of the polynomial
\[
p_{G}\left(x\right)=x^{r}+\alpha_{r-1}x^{r-1}+\cdots+\alpha_{1}x+\alpha_{0},
\]
where the coefficients $\alpha_{0},\ldots,\alpha_{r-1}$ satisfy
\[
\left[\!\begin{array}{c}
\alpha_{0}\\
\alpha_{1}\\
\vdots\\
\alpha_{r-1}
\end{array}\!\right]\!=\!-\!\left[\begin{array}{cccc}
\eigsum_0 & \eigsum_1 & \cdots & \eigsum_{r-1}\\
\eigsum_1 & \eigsum_2 & \cdots & \eigsum_r\\
\vdots & \vdots & \ddots & \vdots\\
\eigsum_{r-1} & \eigsum_r & \cdots & \eigsum_{2r-2}
\end{array}\right]^{-1}\!\left[\begin{array}{c}
\eigsum_r\\
\eigsum_{r+1}\\
\vdots\\
\eigsum_{2r-1}
\end{array}\right]\!.
\]
\end{thm}

\begin{pf}
See Appendix~\ref{app:thm_multi_id_proof}.
\end{pf}

\begin{rmk}
Theorem~\ref{thm:id_multi} provides a methodology for the reconstruction of the observable spectrum of the unknown graph matrix $G$ from $2n$ output observations. From a computational point of view, this method involves the inversion of a lower triangular $2r\times 2r$ matrix, the inversion of an $r\times r$ Hankel matrix, and finding the roots of a degree-$r$ polynomial, where $r$ is the rank of $H$.
\end{rmk}

\section{Continuous-Time Dynamics}\label{sec:ct}

In the case of continuous-time dynamics, there are some subtle but important differences to that of discrete-time. Fortunately, similar results can still be derived in this domain, as we will describe in the following subsections.

\subsection{Network of Single Integrators}\label{subsec:ct_single}

We begin our exposition by considering the case of a network of coupled continuous-time single integrators:
\begin{align}
\begin{aligned}\label{eq:CT_system}
\dot{\mathbf{x}}(t) & = G\mathbf{x}(t),\;\mathbf{x}\left(0\right)=\mathbf{x}_0,\\
y(t) & =\mathbf{c}^{\intercal}\mathbf{x}(t),
\end{aligned}
\end{align}
where $G$ is a graph matrix whose connectivity structure matches that of a potentially weighted and/or directed graph $\G$. We thus have $y\left(t\right)=\mathbf{c}^{\intercal}e^{Gt}\mathbf{x}_0$. In practice, we consider discrete samples $y_k$ of the output with an arbitrary period $\tau>0$, i.e., $y_k \coloneqq y(k\tau)$ for $k \in \mathbb{N}$. Using the Jordan decomposition $G = VJV^{-1}$, we have
\begin{align*}
y_{k} &= \mathbf{c}^{\intercal}V e^{J k\tau}V^{-1}\mathbf{x}_0.
\end{align*}
In contrast to the discrete-time case, here the observations are comprised of exponentiated Jordan matrices, where the $m_i \times m_i$ exponentiated Jordan block $e^{Jk\tau}$ is the upper-triangular matrix
\begin{align}\label{eq:Jordan_exp_ct}
e^{Jk\tau} \!=\! \left[\!\begin{array}{ccccc}
\!e^{\lambda_{i}k\tau} & k\tau e^{\lambda_{i}k\tau} & \!\cdots\! & \frac{(k\tau)^{(m_i-1)}}{(m_i-1)!}e^{\lambda_{i}k\tau}\! \\
 & e^{\lambda_{i}k\tau} & \!\cdots\! & \frac{(k\tau)^{(m_i-2)}}{(m_i-2)!}e^{\lambda_{i}k\tau}\! \\
 &  & \ddots &  \vdots\\
 &  &  & e^{\lambda_{i}k\tau}
\end{array}\!\right]\!.
\end{align}
Thus, our discrete observations may be expressed as
\begin{align}\label{eq:obs_ct}
y_{k} = \sum_{i=1}^{d}\sum_{s=0}^{{m}_i-1}\omega_{i}^{(s)} \frac{(k\tau)^{s}}{s!}(e^{\lambda_i \tau})^k,
\end{align}
with $\omega_{i}^{(s)}$ as defined in \eqref{eq:weights}. Similarly to the discrete-time case, the set of observable eigenvalues is $S_{G} \coloneqq \left\{ \lambda_i \in \sigma(G) \colon \exists s \text{ s.t. } \bar{\omega}_i^{(s)} \neq 0 \right\} $ which, with an arbitrary random initial condition $\x_0$, is almost surely the set of observable eigenmodes of the pair $(G,\mathbf{c}^\intercal)$ according to the PBH test. Fortunately, we may apply analogous results to those in Section~\ref{sec:dt} in order to estimate the eigenvalues corresponding to observable eigenmodes. This notion is formalized in the corollary below.

\begin{cor}\label{cor:CT_single}
Given the sequence of observations $(y_k)_{k=0}^{2n-1}$ from the continuous-time system in~\eqref{eq:CT_system} with fixed sampling rate $\tau>0$, the observable eigenvalues of the graph matrix $G$  may be obtained via
\[
 \lambda_i = \log(\eta_i) / \tau,
\]
where $\eta_i$ are the roots of the polynomial
\[
p_{G}\left(x\right)=x^{r}+\alpha_{r-1}x^{r-1}+\cdots+\alpha_{1}x+\alpha_{0},
\]
whose coefficients $\alpha_0,\ldots,\alpha_{r-1}$ are obtained from the observations $(y_k)_{k=0}^{2n-1}$ as in Theorem~\ref{thm:id}.
\end{cor}
\begin{pf}
See Appendix~\ref{app:CT_single_pf}.
\end{pf}

\subsection{Network of Identical Continuous-Time Agents}\label{subsec:ct_multi}

Similarly to the setting in Section~\ref{subsec:dt_multi}, we consider the dynamics of a network of continuous-time agents beyond single integrators. Assume that each agent is a linear system with state matrix $A$, whose input is a linear combination of the state of its neighbors, its initial state is proportional to a vector $\boldb$, and the measured output is the linear combination $y(t) = \sum_{i}c_i\boldg^\intercal \x_i(t)$. Hence, the global dynamics of the network can be described (in a compact form) analogously to the discrete-time case as 
\begin{align*}
\dot{\mathbf{x}}(t) & =\left(I_{n}\otimes A + G\otimes I_{d}\right)\mathbf{x}(t),\;\mathbf{x}\left(0\right)=\mathbf{x}_0\otimes\boldb,\\
y(t) & =\left(\mathbf{c}\otimes\boldg\right)^{\intercal}\mathbf{x}(t).
\end{align*}
Hence, considering a sampling period $\tau>0$, we have
\begin{align*}
y_k \coloneqq y\left(k\tau\right) & =\left(\mathbf{c}\otimes\boldg\right)^{\intercal}e^{\left(I_{n}\otimes A + G\otimes I_{d}\right)k\tau}\left(\mathbf{x}_0\otimes\boldb\right)\\
 & =\left(\mathbf{c}\otimes\boldg\right)^{\intercal}\left(e^{Gk\tau}\otimes e^{Ak\tau}\right)\left(\mathbf{x}_0\otimes\boldb\right),
\end{align*}
Where the last equality follows by commutativity of the identity matrix and properties of the Kronecker product~\citep[see][]{Petersen12}. Thus,
\begin{align*}
y_k & =\left(\mathbf{c}^{\intercal}V e^{J k\tau}V^{-1}\mathbf{x}_0\right)\left(\boldg^{\intercal}e^{Ak\tau}\boldb\right)\\
 & =\nu_{k}\sum_{i=1}^{d}\sum_{s=0}^{{m}_i-1}\omega_{i}^{(s)} \frac{(k\tau)^{s}}{s!}(e^{\lambda_i \tau})^k,
\end{align*}
where $\nu_{k}\coloneqq\boldg^{\intercal}e^{Ak\tau}\boldb$ and $\omega_{i}^{(s)}$ is defined in \eqref{eq:weights}. Combining Theorem~\ref{thm:id_multi} and Corollary~\ref{cor:CT_single} allows us to use the output measurements $(y_k)_{k=0}^{2n-1}$ to find the roots of a polynomial, which correspond to the values $\eta_i = e^{\lambda_i\tau}$, and we may thus obtain the observable eigenvalues of the graph matrix $G$ after applying a logarithmic transformation.

\section{Simulations}\label{sec:sims}

In this section we illustrate our results, in both discrete- and continuous-time, on networks where the underlying network structure is unknown to us. The evolution of the dynamics of these systems are simulated with an arbitrary random initial condition vector $\x_0$ and an observability vector $\c$. Both $\x_0$ and $\c$ are unknown to the algorithm. Then, we apply Theorem~\ref{thm:id} to estimate the eigenvalues of $G$ from the sequence of observations $\left(y[k]\right)_{k=0}^{2n-1}$ and compare our estimated eigenvalues against the true spectrum of the graph matrix $G$.

\begin{figure*}[!ht]
\centering
    \begin{subfigure}{.3\textwidth}
    \centering
    \includegraphics[width=0.9\linewidth]{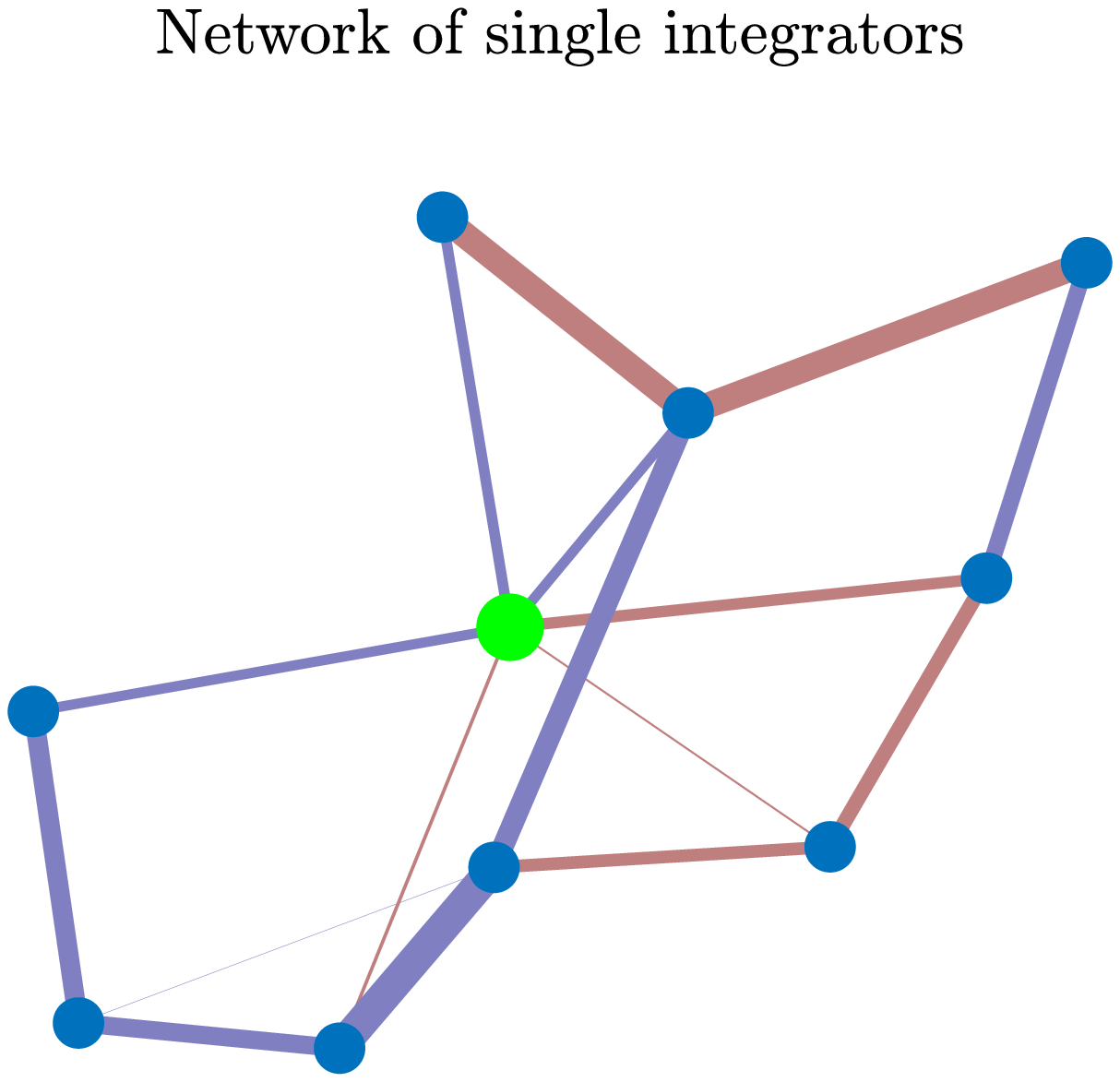}
    \caption{Network topology, with single output agent highlighted. Edge thickness corresponds to edge weight; red edges have negative weights.}
    \label{fig:dt_ntwk}
    \end{subfigure}
    \quad
    \begin{subfigure}{.3\textwidth}
    \centering
    \includegraphics[width=\linewidth]{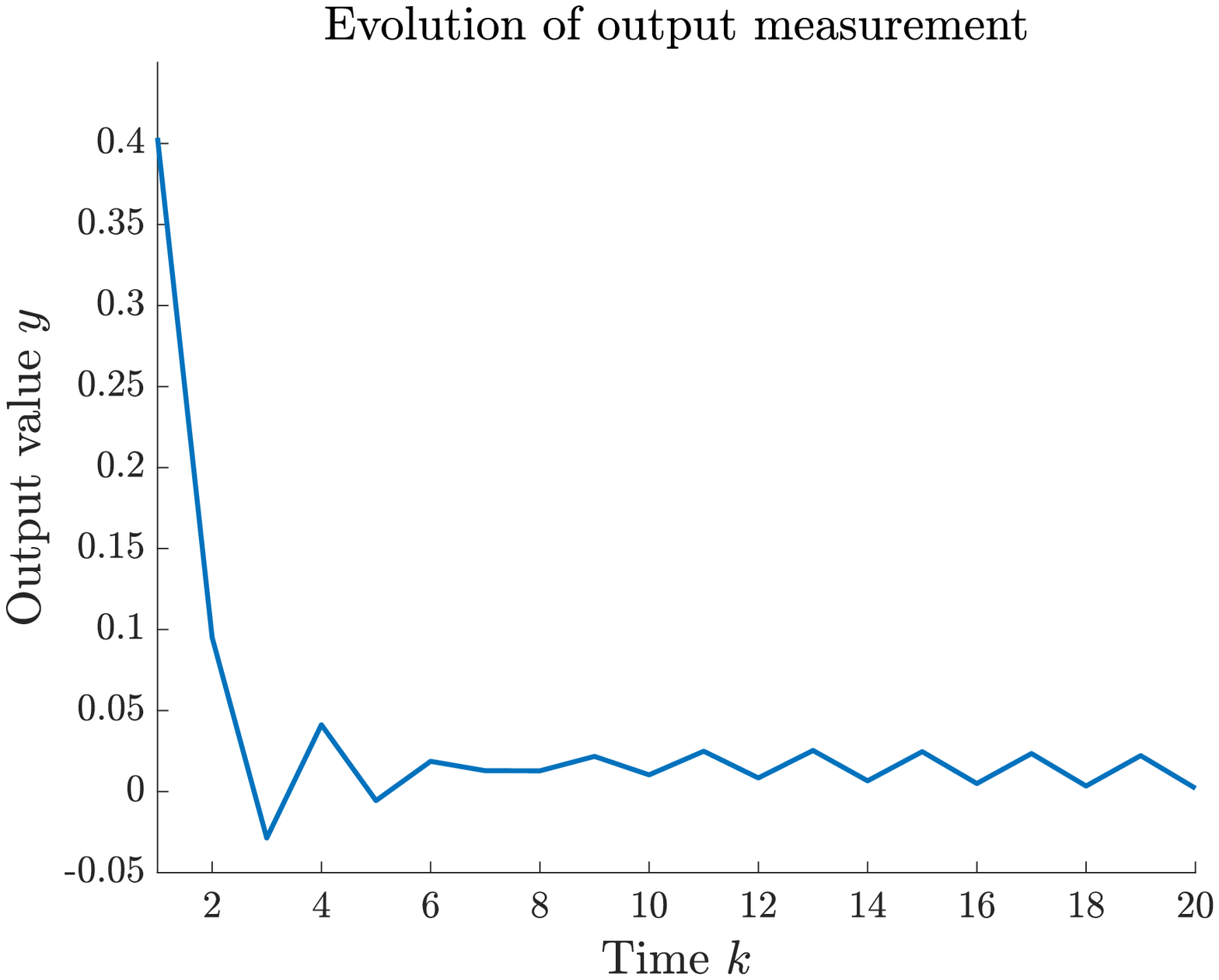}
    \caption{Output $y[k] = e_i^\intercal G^k x_0$, where we observe only agent $i$.\\}
    \label{fig:dt_y}
    \end{subfigure}
    \quad
    \begin{subfigure}{.3\textwidth}
    \centering
    \includegraphics[width=\linewidth]{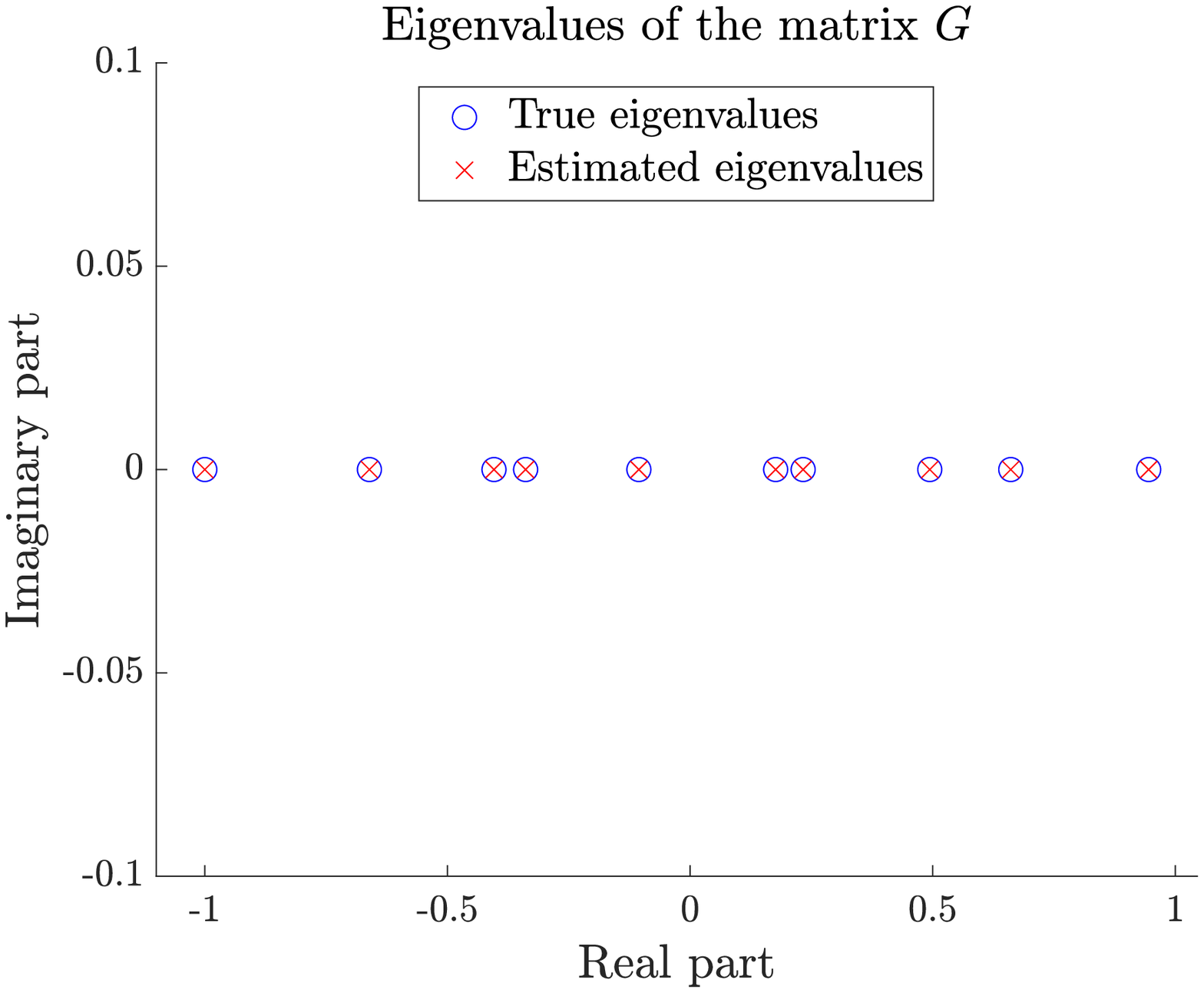}
    \caption{Comparison of true and estimated eigenvalues; repeated values overlaid.\\}
    \label{fig:dt_eigs}
    \end{subfigure}
    \caption{$10$-agent preferential attachment network in discrete-time, generated according to \cite{Barabasi99}. The initial condition is randomly generated as $\x_0 \sim \text{Uniform}[0,1]^n$. There are $10$ eigenvalues of $G$ in this case, which are all recovered via our estimation approach.}
    \label{fig:dt_sims}
\end{figure*}

Figure~\ref{fig:dt_sims} shows the result of using Theorem~\ref{thm:id} on the undirected, randomly generated $10$-agent preferential attachment network shown in Figure~\ref{fig:dt_ntwk} \citep[see][]{Barabasi99}. We model each agent using a single integrator dynamics in discrete-time, as in~\eqref{eq:dt_integrator}. We assume that we only have access to the output of the integrator agent indicated in green in Fig~\ref{fig:dt_ntwk}. The thickness of edges in Figure~\ref{fig:dt_ntwk} is proportional to their weight in the graph matrix $G$, with negative weights shown in red; the weights are generated according to a $\text{Uniform}[-1,1]$ distribution. In Figure~\ref{fig:dt_y}, we show the evolution of the output signal; as only one agent's output is measured, this may be viewed as a decentralized eigenvalue estimation problem. Figure~\ref{fig:dt_eigs} compares both the true and estimated eigenvalues of $G$. In this case there are $10$  eigenvalues of $G$, and all of these are perfectly recovered using a sequence of 20 measurements retrieved from a single agent.

\begin{figure*}[!ht]
\centering
    \begin{subfigure}{.3\textwidth}
    \centering
    \includegraphics[width=0.9\linewidth]{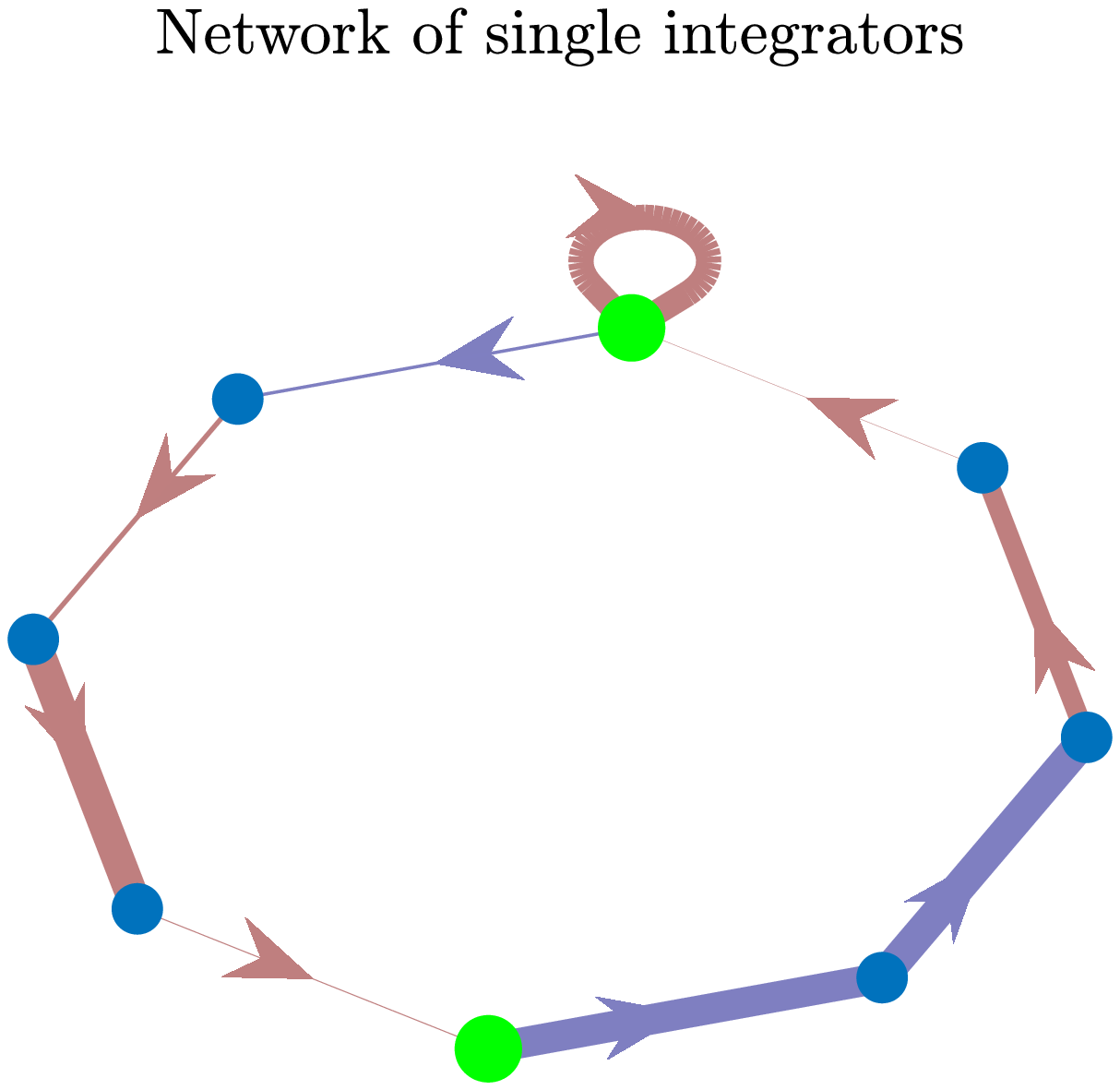}
    \caption{Network topology, with output agents highlighted. Edge thickness corresponds to edge weight; red edges have negative weights.}
    \label{fig:ct_ntwk}
    \end{subfigure}
    \quad
    \begin{subfigure}{.3\textwidth}
    \centering
    \includegraphics[width=\linewidth]{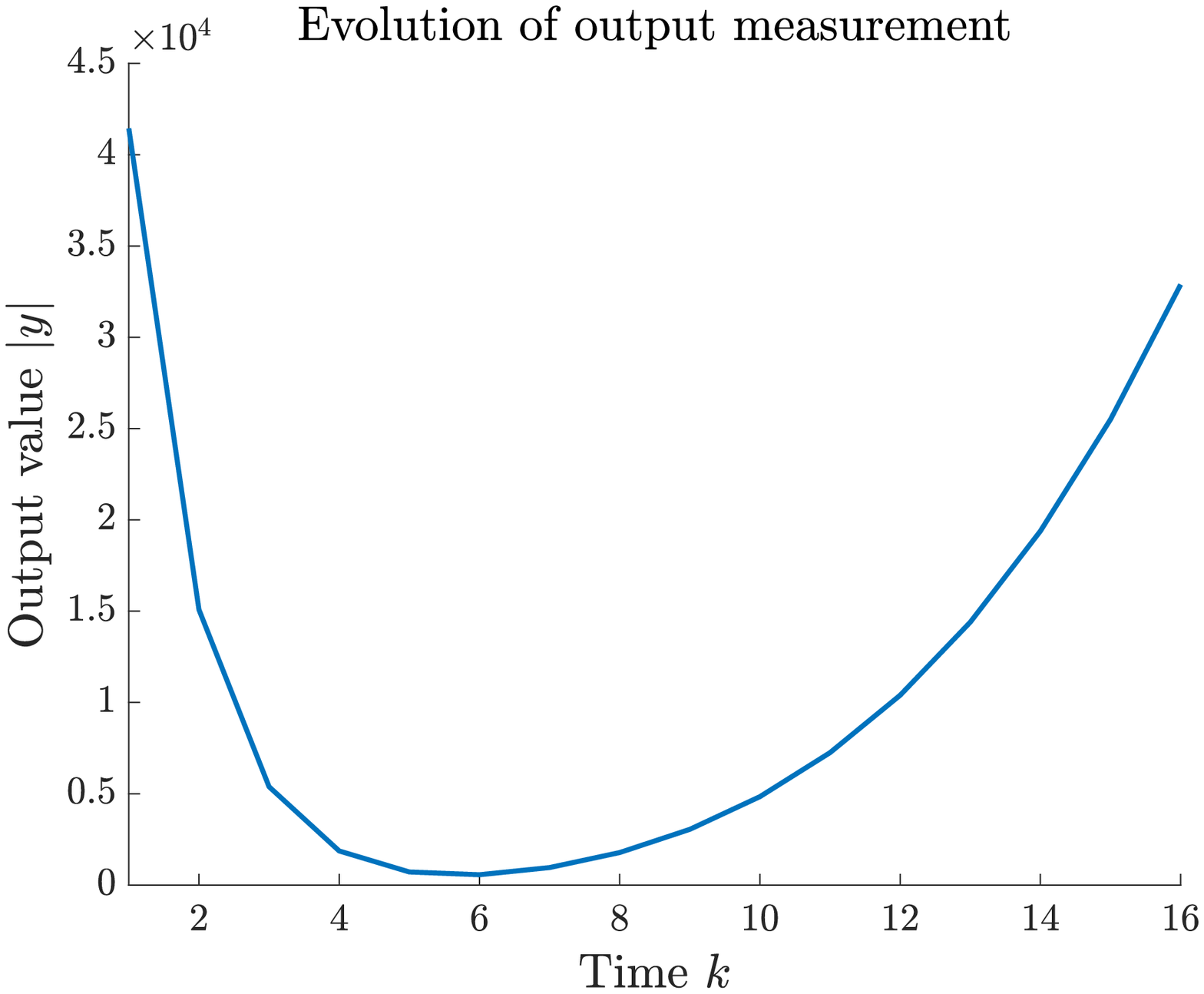}
    \caption{Output $y_k = c^\intercal e^{-Gk\tau} x_0$; agents are observed with equal weight.\\}
    \label{fig:ct_y}
    \end{subfigure}
    \quad
    \begin{subfigure}{.3\textwidth}
    \centering
    \includegraphics[width=\linewidth]{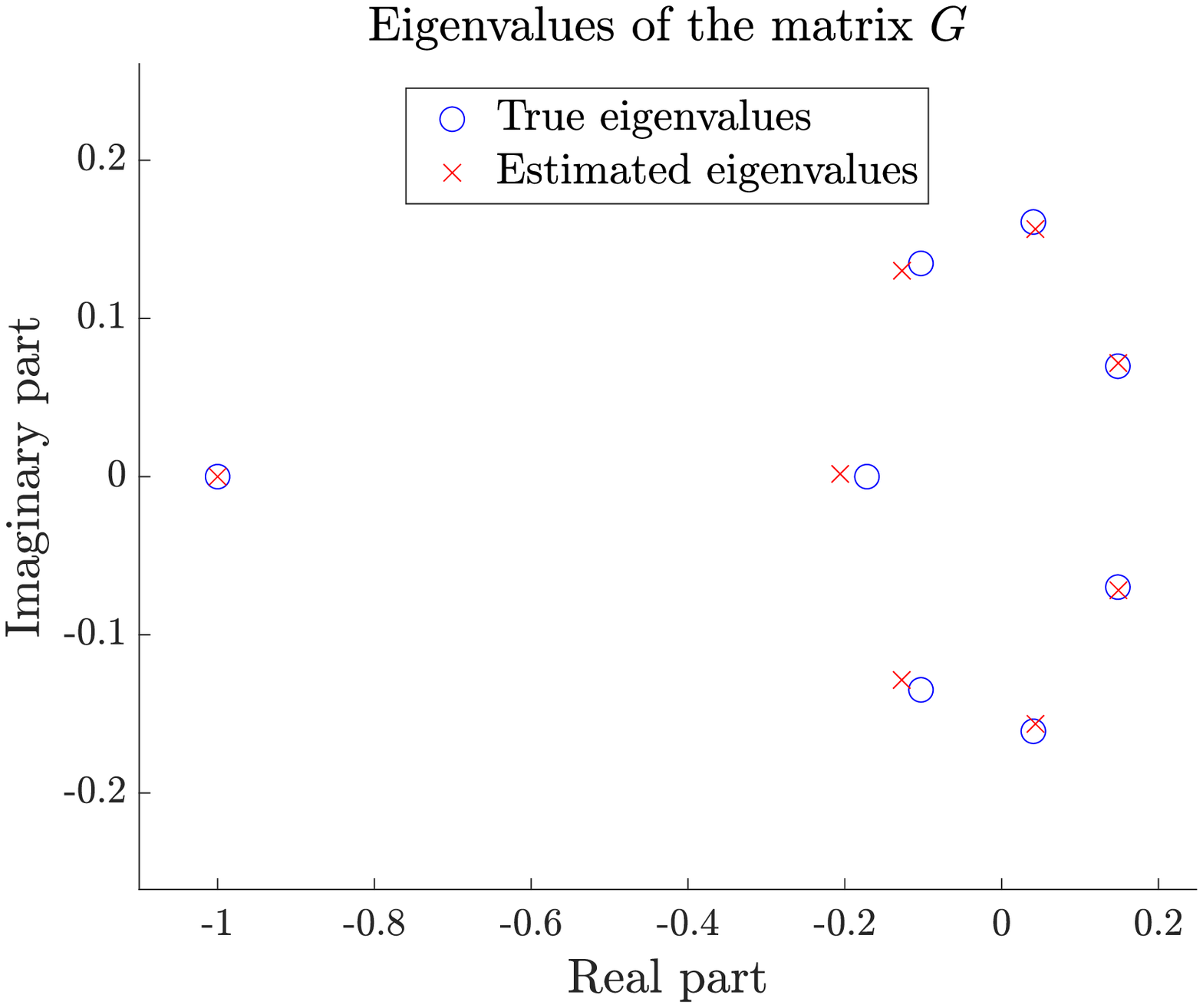}
    \caption{Comparison of true and estimated eigenvalues.\\}
    \label{fig:ct_eigs}
    \end{subfigure}
    \caption{$8$-agent single integrator ring network in continuous-time, with sampling rate $\tau = 1$ and random initial condition $\x_0 \sim \text{Uniform}[0,1]^n$. Here there are $8$ eigenvalues of $G$, all of which are recovered via our estimation approach.}
    \label{fig:ct_sims}
\end{figure*}

In Figure~\ref{fig:ct_sims} we apply our estimation approach on the $8$-agent weighted and directed ring network shown in Figure~\ref{fig:ct_ntwk}, wherein the agents obey the continuous-time dynamics described in Section~\ref{sec:ct}. Again, edge thickness in Figure~\ref{fig:ct_ntwk} corresponds to the edges' weights in the graph matrix $G$, with negative weights shown in red; the weights are generated according to a $\text{Uniform}[-1,1]$ distribution. In this case the output is a linear combination of the states of the two agents highlighted in Figure~\ref{fig:ct_ntwk}. Although our realization of $G$ renders an unstable system and the output eventually grows exponentially (as shown in Figure~\ref{fig:ct_y}), we are still able to recover the entirety of the true spectrum of $G$ with high accuracy as shown in Figure~\ref{fig:ct_eigs}. The difference in accuracy from Figure~\ref{fig:dt_sims} is due to the numerical sensitivity of root-finding techniques, since the outputs are large due to the system being unstable.

\begin{figure*}[!ht]
\centering
    \begin{subfigure}{.3\textwidth}
    \centering
    \includegraphics[width=0.9\linewidth]{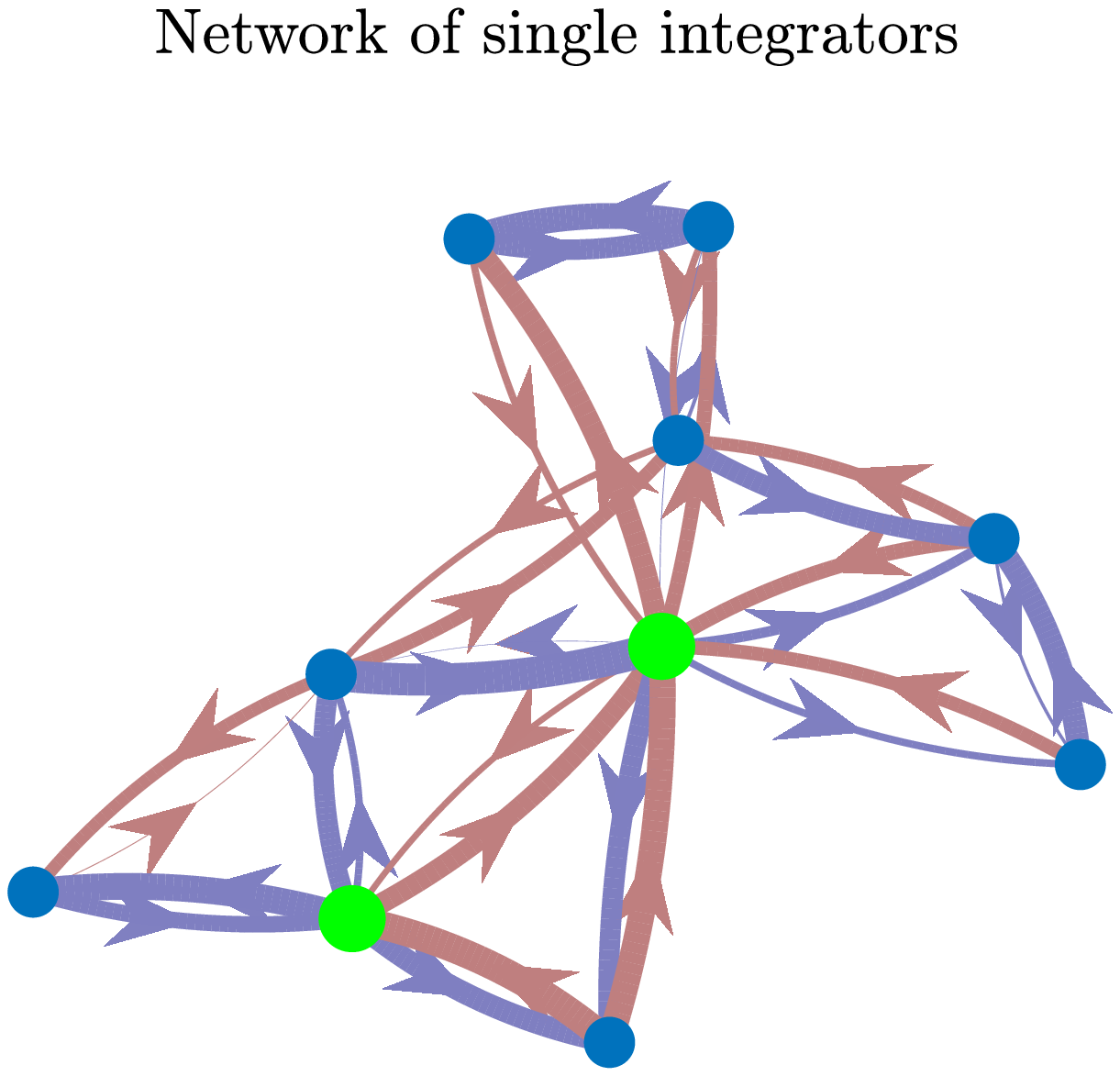}
    \caption{Network topology, with output agents highlighted. Edge thickness corresponds to edge weight; red edges have negative weights.}
    \label{fig:multi_ntwk}
    \end{subfigure}
    \quad
    \begin{subfigure}{.3\textwidth}
    \centering
    \includegraphics[width=\linewidth]{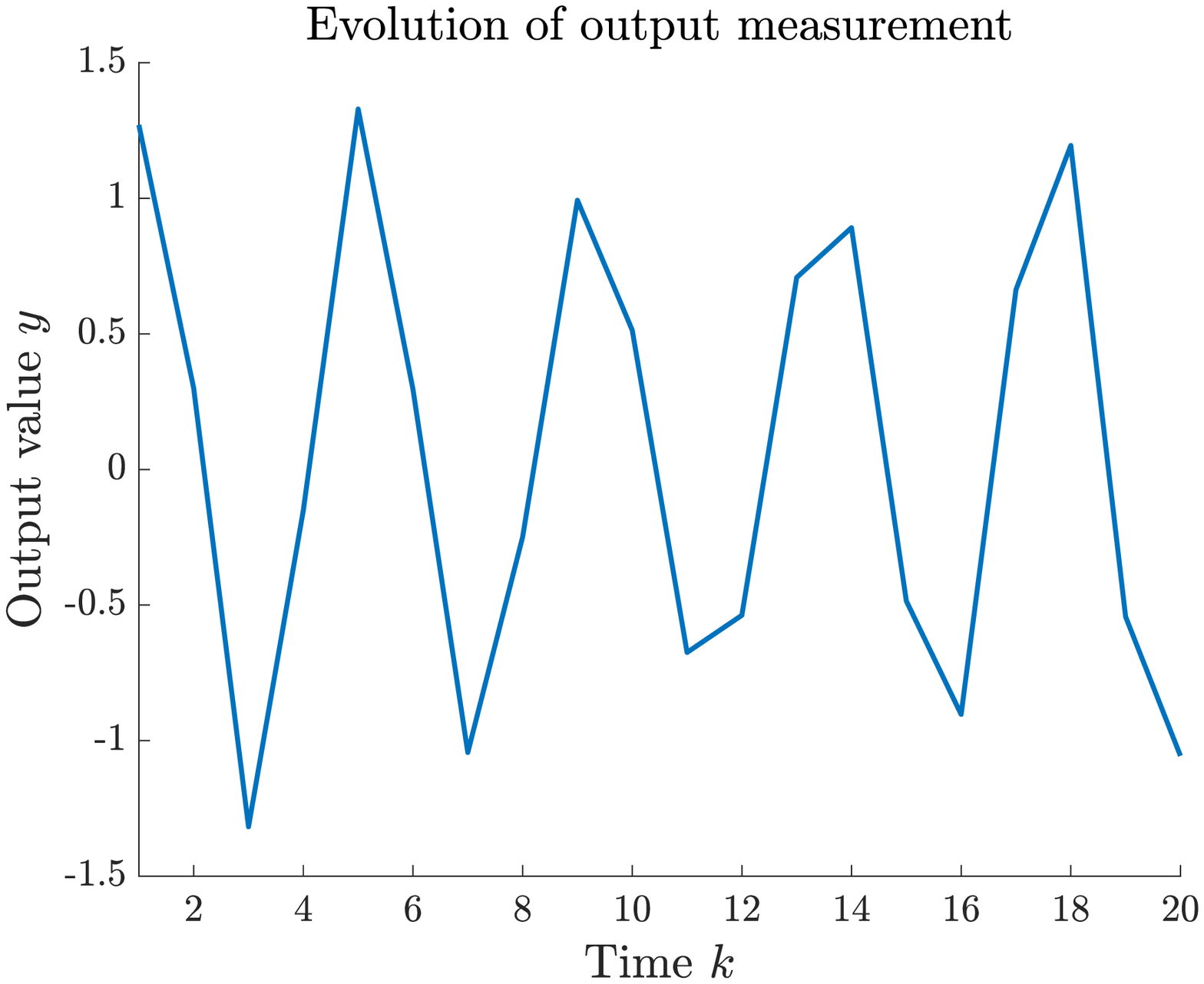}
    \caption{Output $y[k]$ in this case is equal to $\left(\mathbf{c}\!\otimes\!\boldg\right)^{\intercal}\!\left(I_{n}\!\otimes\! A \!+\! G\!\otimes\! I_{d}\right)^k\!(\x_0\!\otimes\boldb)$; agents are observed with equal weight.}
    \label{fig:multi_y}
    \end{subfigure}
    \quad
    \begin{subfigure}{.3\textwidth}
    \centering
    \includegraphics[width=\linewidth]{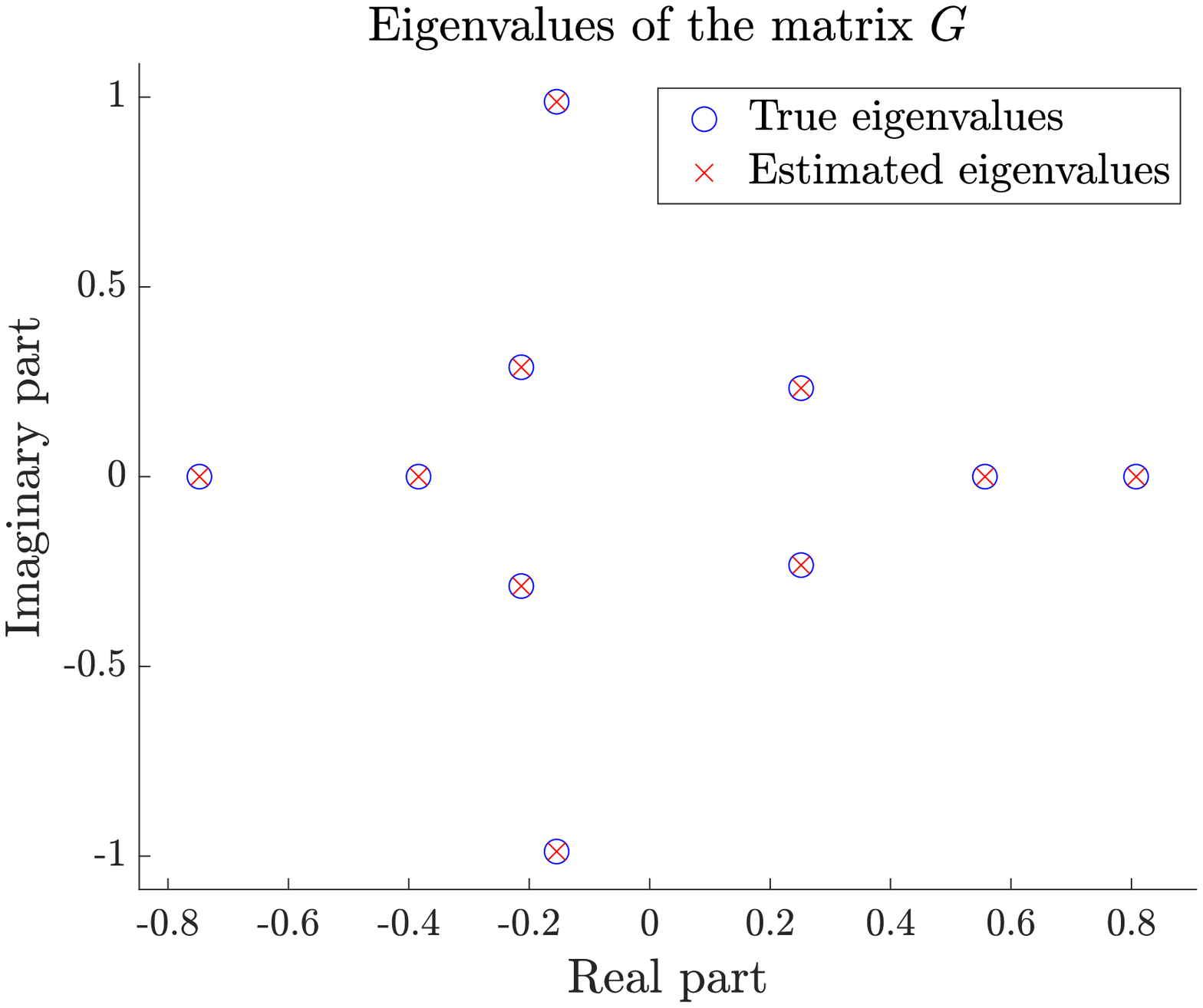}
    \caption{Comparison of true and estimated eigenvalues. \\}
    \label{fig:multi_eigs}
    \end{subfigure}
    \caption{$10$-agent preferential attachment network in discrete-time, generated according to \cite{Barabasi99}. The dynamics here follow the more general case of \eqref{eq:dt_multi} from Section~\ref{subsec:dt_multi}, where each node has a $3$-dimensional state. The initial condition is generated according to $\x_0 \sim \text{Uniform}[0,1]^n$, the vectors $\boldb$ and $\boldg$ are generated according to $\text{Uniform}[0,1]^3$, and the entries of the symmetric matrix $A$ are generated as $a_{ij} \sim \text{Uniform}[0,1]^n,~i\geq j$. There are $10$ eigenvalues of $G$ in this case, which are all recovered via our estimation approach.}
    \label{fig:multi_sims}
\end{figure*}

In Figure~\ref{fig:multi_sims} we estimate the eigenvalues of a network of $10$ discrete-time identical agents. Figure~\ref{fig:multi_ntwk} displays a randomly generated preferential attachment network over which the agents interact. In this case the edges are weighted according to a $\text{Uniform}[-1,1]$ distribution, with thickness representing weight and negatively-weighted edges shown in red. The measurements we observe in Figure~\ref{fig:multi_y} are the sum of the outputs of the two agents located on the green nodes of the network. After applying Theorem~\ref{thm:id_multi}, all $10$ eigenvalues are recovered as shown in Figure~\ref{fig:multi_eigs}.


\section{Conclusion}\label{sec:conc}

In this paper, we have proposed an efficient methodology for estimating the eigenvalues of any arbitrary graph matrix of a network of interacting dynamical agents using a set of dynamical measurements. This graph matrix may be directed, may have edge weights, incorporate self-loops, and may render the system unstable. Unlike other methods, we require only a single finite sequence of discrete, temporal measurements from the multiagent network of length, at most, $2n$. Moreover, we need no prior knowledge of the network topology, initial condition, or which agents are contributing to the measurements. For any arbitrary random initial condition our approach is able to recover all eigenvalues corresponding to observable eigenmodes of the pair $(G,\c^\intercal)$, almost surely. We develop our technique for systems in both discrete- and continuous-time, and consider the case of agents modeled by single integrators as well as more complex dynamics. Our simulation results show that we are able to recover the observable spectrum of the graph matrix in all cases with high accuracy.

\appendix
\section{Proofs}

\subsection{Proof of Lemma~\ref{lem:rank}}\label{app:lem_rank_proof}

Recall the set of indices corresponding to observable eigenvalues $\I = \left\{ i\in\{1,\ldots,n\} : \lambda_i \in \S_G \right\}$, and the total weights corresponding to each unique eigenvalue $\bar{\omega}_{i}^{(s)} = \sum_{j:\lambda_{j}=\lambda_{i}}\omega_{j}^{(s)}$ from~\eqref{eq:total_weights}. Now let $\mathbf{v}_{i}\coloneqq\left[1,\lambda_{i},\lambda_{i}^{2}\ldots,\lambda_{i}^{n-1}\right]$ and $\bin{n}{k} \coloneqq \binom{n}{k}$. Then combining~\eqref{eq:Jordan_exp}, \eqref{eq:obs_total}, and~\eqref{eq:hankel} we obtain
\begin{align*}
H \!&= \!\sum_{i\in\I}\sum_{s=0}^{\hat{m}_i-1} \bar{\omega}_i^{(s)} \!\left[\!\begin{array}{cccc}
\bin{0}{s} & \bin{1}{s}\lambda_i^{1-s} & \cdots & \bin{n-1}{s}\lambda_i^{n-1-s} \\
\bin{1}{s}\lambda_i^{1-s} & \bin{2}{s}\lambda_i^{2-s} & \cdots & \bin{n}{s}\lambda_i^{n-s} \\
 \vdots & \vdots & \ddots & \vdots \\
\bin{n-1}{s}\lambda_i^{n-1-s} & \bin{n}{s}\lambda_i^{n-s} & \cdots & \bin{2n-2}{s}\lambda_i^{2n-2-s} \\
\end{array}\!\right] \\
    &=\sum_{i\in\I}\sum_{s=0}^{\hat{m}_i-1} \frac{\bar{\omega}_{i}^{(s)}}{s!} \frac{d^s}{d\lambda_i^s} \left(\mathbf{v}_{i}\mathbf{v}_{i}^{\intercal}\right) =: \sum_{i\in\I} H_i,
\end{align*}
where the derivative is taken element-wise to the entries of the matrix $\mathbf{v}_{i}\mathbf{v}_{i}^{\intercal}$. Notice that for all $s$ and any given $i,j\in\I$ the Hankel matrices $\frac{d^s}{d\lambda_i^s} \left(\mathbf{v}_{i}\mathbf{v}_{i}^{\intercal}\right)$ and $\frac{d^s}{d\lambda_j^s} \left(\mathbf{v}_{j}\mathbf{v}_{j}^{\intercal}\right)$ have orthogonal ranges since the $\lambda_i$ for $i \in \I$ are unique, and so $\mathbf{v}_i$ and $\mathbf{v}_j$ are linearly independent.

Let us now examine the ranks of the matrices $D_i^{(s)} \coloneqq \frac{d^s}{d\lambda_i^s} \left(\mathbf{v}_{i}\mathbf{v}_{i}^{\intercal}\right)$ for a particular $i \in \I$. We will proceed via induction on $s$ to show that $\text{rk}(D_i^{(s)}) = s+1$. For the base case of $s=0$ we have $D_i^{(0)} = \mathbf{v}_i\mathbf{v}_i^\intercal$, which clearly has rank $1$.

Now we assume $\text{rk}(D_i^{(s-1)}) = s$. The $j$-th column of $D_i^{(s)}$ is of the form
\[
d_{i,j}^{(s)} \coloneqq s!
\begin{bmatrix}
\bin{j-1}{s}\lambda^{j-1-s} \\
\bin{j}{s}\lambda^{j-s} \\
\vdots \\
\bin{j + n - 2}{s}\lambda^{j + n - 2 - s}
\end{bmatrix}
=
s!\lambda^{j-1-s}
\begin{bmatrix}
\bin{j-1}{s} \\
\bin{j}{s}\lambda \\
\vdots \\
\bin{j + n - 2}{s}\lambda^{n-1}
\end{bmatrix}\!.
\]
Recall that $\bin{k}{s} = \binom{k}{s} = 0$ for $ k < s$. By the leading-zero structure of $D_i^{(s)}$, wherein the first column has $s$ leading zeros followed by a nonzero value, the second has $s-1$ leading zeros followed by a nonzero value, all the way to the $s$-th column having a nonzero value in the first component, we can see that $\text{rk}(D_i^{(s)}) \geq s+1$. Now take any collection of $s+2$ columns of $D_i^{(s)}$, and we will show they must be linearly dependent. Via the identity
\[
\binom{j}{s} - \binom{j-k}{s} = \binom{j-1}{s-1} + \binom{j-2}{s-1} + \cdots + \binom{j-k}{s-1},
\]
we may write
\begin{align*}
\frac{d_{i,j}^{(s)}}{\lambda_i^{j-1-s}} - \frac{d_{i,j-k}^{(s)}}{\lambda_i^{j-k-1-s}} = \frac{s!}{(s+1)!}\sum_{l=1}^k\frac{d_{i,j-l}^{(s-1)}}{\lambda_i^{j-l-1-s}}.
\end{align*}
In other words, we may express the $j$-th and $(j-k)$-th columns of $D_i^{(s)}$ as a linear combination of exactly $k$ columns from $D_i^{(s-1)}$. Since we have a collection of $s+2$ columns of $D_i^{(s)}$, we will need at least $s+1$ unique columns of $D_i^{(s-1)}$ to express linear combinations of our entire collection (in the case where the columns are sequential), but may need more. However, the rank of  $D_i^{(s-1)}$ is $s$, so any collection of at least $s+1$ unique columns of $D_i^{(s-1)}$ must be linearly dependent; hence, our collection of $s+2$ columns of $D_i^{(s)}$ must be linearly dependent. Thus, $\text{rk}\left(D_i^{(s)}\right) = s+1$.

We will now examine the ranges of the matrices $D_i^{(s)}$ for a particular $i \in \I$. For $1 \leq s < j-k$ and $1 \leq k < j$ we have the identity 
\begin{align*}
\binom{j}{s} - \binom{j-k}{s} = \binom{j-k}{s-k}.
\end{align*}
Thus,
\[
d_{i,j-k}^{(s-k)} = \frac{(s-k)!}{s!}\left[d_{i,j}^{(s)} - \lambda_i^kd_{i,j-k}^{(s)}\right]\!.
\]
In other words, we may write the $(j-k)$-th column of $D_i^{(s-k)}$ as a linear combination of the $j$-th and $(j-k)$-th columns of $D_i^{(s)}$ for $1 \leq k \leq s$ and $k < j \leq n$. Recall that $\text{rk}(D_i^{(s-k)}) = s-k+1$. Since we may write the first $s - k + 1$ columns of $D_i^{(s-k)}$ as linear combinations of the columns of $D_i^{(s)}$, the same is true for all columns of $D_i^{(s-k)}$. Thus $\text{range}(D_i^{(s-k)}) \subseteq \text{range}(D_i^{(s)})$ for $1 \leq k \leq s \leq \tilde{m}_i - 1$. Hence,
\begin{align*}
	&\text{range}\left(H_i\right) = \text{range}\!\left(\sum_{s=0}^{\hat{m}_i-1}\frac{\bar{\omega}_{i}^{(s)}}{s!}D_i^{(s)}\!\right) \!= \text{range}(D_{i}^{(\tilde{m}_i - 1)}) \\
	&\Rightarrow\quad \text{rk}(H_i) = \text{rk}(D_i^{(\tilde{m}_i-1)}),
\end{align*}
where $\tilde{m}_i$, as defined in~\eqref{eq:m_tilde}, is the largest index with a nonzero total weight. Thus, the rank of $H_i$ is simply the largest $s$ for which $\bar{\omega}_i^{(s-1)} \neq 0$, i.e., $\text{rk}(H_i) = \tilde{m}_i$.

Since for all $s$ and any $i \neq j \in \I$ the matrices $D_i^{(s)}$ and $D_j^{(s)}$ have orthogonal ranges, we have that $\text{rg}(H) = \text{rg}(\sum_{i \in \I} H_i) = \oplus_{i\in\I}\text{rg}(H_i)$, and hence $\text{rk}(H) = \sum_{i \in \I}\text{rk}(H_i)$. Therefore, the rank of $H$ is equal to the sum of the sizes of the largest observable Jordan blocks for each unique eigenvalue, which is
$\sum_{i \in \I} \tilde{m}_i$. \qed

\subsection{Proof of Theorem~\ref{thm:id}}\label{app:thm_proof}

By definition, we know that at most we may recover all eigenvalues corresponding to observable eigenmodes, i.e., $\lambda_i\in\S_{G}$. As before, let $\I = \left\{ i\in\{1,\ldots,n\} : \lambda_i \in \S_{G}\right\}$. By Lemma \ref{lem:rank}, we know that $\text{rk}\left(H\right) = \sum_{i \in \I} \tilde{m}_i$, which we denote by $r$. Define the following polynomial:
\[
p_{G}\left(x\right) \coloneqq \prod_{i\in\I}(x-\lambda_i)^{\tilde{m}_i} =  x^{r}+\alpha_{r-1}x^{r-1}+\cdots+\alpha_{1}x+\alpha_{0},
\]
where $\tilde{m}_i$ is defined in~\eqref{eq:m_tilde}. Notice that, since the eigenvalues are unknown, the coefficients of the polynomial are also unknown. In what follows, we propose an efficient technique to find these coefficients. 

Let us calculate $p_{G}(J_i)$ for each $i \in \I$. Recall that there may be multiple Jordan blocks associated with a single eigenvalue, and that the Jordan block $J_l$ is of size $m_l\times m_l$. First consider the case that there exists some Jordan block $i$ such that $m_i = \tilde{m}_i$.
By Cayley-Hamilton theorem we know $(J_i - \lambda_i I_{m_i})^{\tilde{m}_i} = \mathbf{0}_{m_i\times m_i}$, and so\begin{align*}
	p_G(J_i) = J_i^r + \alpha_{r-1}J_i^{r-1} + \cdots + \alpha_1 J_1 + \alpha_0 = \mathbf{0}_{m_i\times m_i}.
\end{align*}
Note from~\eqref{eq:Jordan_exp} that each upper diagonal of $J_i^r$ contains the same values. For ease of exposition, define $\bin{n}{k} = \binom{n}{k}$. Hence, since $J_i^r$ is of size $m_i \times m_i$, we in fact have $m_i$ separate equations (one per upper diagonal) of the form
\[
	\bin{r}{s}\lambda_i^{r\!-\!s} \!+\! \alpha_{r\!-\!1}b^{r\!-\!1}_{s}\lambda_i^{(r\!-\!s) - 1} \!+\! \cdots \!+\! \alpha_{s+1}\bin{s+1}{s}\lambda_i + \alpha_{s} \!=\! 0,
\]
for $s \in \{0,\ldots,m_i-1\}$. If there is no Jordan block $i$ such that $m_i = \tilde{m}_i$, then pick one such that $m_i > \tilde{m}_i$, and consider the first $\tilde{m}_i$ upper diagonals of $(J_i - \lambda_i I_{m_i})^{\tilde{m}_i}$, which will be zero. Multiplying the equations above by the corresponding total weights $\bar{\omega}_i^{(s)}$, some of which may be zero, we obtain for $s \in \{0,\ldots,\tilde{m}_i-1\}$
\[
	\bar{\omega}_i^{(s)} \!\left( \bin{r}{s}\lambda_i^{r\!-\!s} \!+\! \alpha_{r\!-\!1}\bin{r\!-\!1}{s}\!\lambda_i^{(r\!-\!s) - 1} \!+\! \cdots \!+\! \alpha_{s+1}\bin{s+1}{s}\!\lambda_i \!+\! \alpha_{s} \right) \!=\! 0.
\]
Summing all of these equations, noting that $\bin{r}{s} = \binom{r}{s} = 0$ for $r < s$, defining $\alpha_r = 1$, we have
\begin{align*}
	&\!\sum_{s=0}^{\tilde{m}_i-1}\!\bar{\omega}_i^{(s)} \!\left( \bin{r}{s}\lambda_i^{r\!-\!s} \!+\! \alpha_{r\!-\!1}\bin{r\!-\!1}{s}\lambda_i^{(r\!-\!s) - 1} \!+\! \cdots \!+\! \alpha_{s+1}\bin{s+1}{s}\lambda_i \!+\! \alpha_{s} \right)  \\
	&= \sum_{s=0}^{\tilde{m}_i-1}\bar{\omega}_i^{(s)} \sum_{l=0}^r \alpha_l \binom{l}{s} \lambda_i^{l-s} = 0.
\end{align*}
Now, let us sum over all eigenvalues $\lambda_i \in \S_G$:
\begin{align*}
&\sum_{i\in\I} \sum_{s=0}^{\tilde{m}_i-1}\bar{\omega}_i^{(s)} \sum_{l=0}^r \alpha_l \binom{l}{s} \lambda_i^{l-s} \\
&= \sum_{l=0}^r \alpha_l \sum_{i\in\I}\sum_{s=0}^{\tilde{m}_i-1}\bar{\omega}_i^{(s)} \binom{l}{s}  \lambda_i^{l-s} = \sum_{l=0}^r \alpha_l y_l = 0,
\end{align*}
by definition of the observations $y_s$ from~\eqref{eq:obs}. Now, for $k \in \{1,\ldots,r-1\}$, let us examine the equations
\[
	J^k p_\G(J) = J^{r+k} + \alpha_{r-1}J^{r+k-1} + \cdots + \alpha_1 J^{k+1} + \alpha_0J^k.
\]
Repeating the same process from above, we obtain for $k \in \{1,\ldots,r-1\}$
\begin{align*}
&\sum_{i \in \I}\sum_{s=0}^{\tilde{m}_i-1}\bar{\omega}_i^{(s)} \sum_{l=0}^r \alpha_l \binom{l+k}{s}  \lambda_i^{l+k-s} \\
&= \sum_{l=0}^r \alpha_l \sum_{i \in \I}\sum_{s=0}^{\tilde{m}_i-1}\bar{\omega}_i^{(s)} \binom{l+k}{s}  \lambda_i^{l+k-s} = \sum_{l=0}^r \alpha_l y_{l+k} = 0.
\end{align*}
In summary, we have $r$ equations of the form
\[
	y_{r+k} + \alpha_{r-1}y_{r+k-1} + \cdots + \alpha_1 y_{k+1} + \alpha_0 y_k = 0,
\]
where $k \in \{0,\ldots,r-1\}$. In matrix form,
\[
\left[\begin{array}{cccc}
y_0 & y_1 & \cdots & y_{r-1} \\
y_1 & y_2 & \cdots & y_r\\
\vdots & \vdots & \ddots & \vdots\\
y_{r-1} & y_r & \cdots & y_{2r-2}
\end{array}\right]\left[\begin{array}{c}
\alpha_{0}\\
\alpha_{1}\\
\vdots\\
\alpha_{r-1}
\end{array}\right]=-\left[\begin{array}{c}
y_r\\
y_{r+1}\\
\vdots\\
y_{2r-1}
\end{array}\right]\!.
\]
By Lemma~\ref{lem:rank} we know $\text{rk}(H) = r$ and hence we may find the values of the coefficients $\alpha_0,\ldots,\alpha_{r-1}$ by a simple matrix inversion. Using these coefficients we can compute the roots of $p_{G}$ to recover the eigenvalues of $G$ that are in the set $\S_{G}$, i.e., those eigenvalues $\lambda_i$ corresponding to the observable eigenmodes of the dynamics. Moreover, the multiplicity of the root $\lambda_i$ will be $\tilde{m}_i$; hence, we recover $\lambda_i$ with multiplicity of exactly $\tilde{m}_i$. \qed

\subsection{Proof of Theorem~\ref{thm:id_multi}}\label{app:thm_multi_id_proof}

Considering the Jordan decomposition $G = V J V^{-1}$, we have
\begin{align*}
&\left(I_{n}\otimes A+G\otimes I_{d}\right)^{k} \\
 &~= \left[\left(V\otimes I_{d}\right)\left(I_{n}\otimes A+ J\otimes I_{d}\right)\left(V^{-1}\otimes I_{d}\right)\right]^k \\
 &~=\left(V\otimes I_{d}\right)\left(I_{n}\otimes A+J\otimes I_{d}\right)^{k}\left(V^{-1}\otimes I_{d}\right)\\
 &~= \left(V\otimes I_{d}\right) \! \left[\sum_{s=0}^{k}\binom{k}{s}\left(I_{n}\otimes A^{k-s}\right)\left(J^{s}\otimes I_{d}\right)\right] \! \left(V^{-1}\otimes I_{d}\right).
\end{align*}
Thus,
\begin{align*}
&y\left[k\right] = \left(\mathbf{c}\otimes\boldg\right)^{\intercal}\left(I_{n}\otimes A+ G\otimes I_{d}\right)^{k}\left(\mathbf{x}_0\otimes\boldb\right) \\
&=\sum_{s=0}^{k}\!\binom{k}{s}\! \left(\mathbf{c}^{\intercal}V\otimes\boldg^{\intercal}\right) \! \left(I_n\otimes A^{k-s}\right) \! \left(J^{s}\otimes I_{d}\right) \! \left(V^{-1}\mathbf{x}_0\otimes\boldb\right)\\
 & =\sum_{s=0}^{k}\binom{k}{s}\left(\mathbf{c}^{\intercal}V J^{s} V^{-1}\mathbf{x}_0\right)\left(\boldg^{\intercal}A^{k-s}\boldb\right).
\end{align*}
Hence, we obtain
\begin{align}\label{eq:obs_multi}
y\left[k\right] &= \sum_{s=0}^{k}\binom{k}{s}\nu_{k-s}\sum_{i=1}^{d}\sum_{s=0}^{m_{i}-1}\omega_{i}^{(s)}{k \choose s}\lambda_{i}^{k-s} \cr
	&= \sum_{s=0}^{k}\binom{k}{s}\nu_{k-s}\eigsum_{s},
\end{align}
where $\eigsum_s = \sum_{i=1}^{d}\sum_{s=0}^{m_{i}-1}\omega_{i}^{(s)}{k \choose s}\lambda_{i}^{k-s}$ and $\nu_{k-s} = \boldg^{\intercal}A^{k-s}\boldb$. From the sequence $\left(y\left[k\right]\right)_{k=0}^{2n-1}$, we obtain a lower triangular system of linear equations that can be solved to find the sequence $\left(m_{k}\right)_{k=0}^{2n-1}$. Specifically, if we collect $2r$ observations, with $\bin{k}{s} = \binom{k}{s}$, we have that \eqref{eq:obs_multi} for $k=0,\ldots,2r-1$ results in
\begin{align*}
\!\begin{bmatrix}
y_{0}\\
y_{1}\\
\vdots\\
\!y_{2r-1}\!
\end{bmatrix} \!=\! \begin{bmatrix}
\bin{0}{0}\nu_{0} & 0 & \cdots & 0\\
\bin{1}{0}\nu_{1} & \bin{1}{1}\nu_{0} & \cdots & 0\\
\vdots & \vdots & \ddots & \vdots\\
\!\bin{2r-1}{0}\nu_{2r-1} & \bin{2r-1}{1}\nu_{2r-2} \!&\! \cdots \!&\! \bin{2r-1}{2r-1}\nu_{0}\!
\end{bmatrix}\!\begin{bmatrix}
\eigsum_{0}\\
\eigsum_{1}\\
\vdots\\
\!\eigsum_{2r-1}\!
\end{bmatrix}
\end{align*}
As long as $\nu_0 = \boldg^\intercal\boldb \neq 0$, the above matrix is full-rank. We may then recover the values $\eigsum_s$ by a simple inversion, and apply Theorem~\ref{thm:id} to find the eigenvalues of $G$. \qed

\subsection{Proof of Corollary~\ref{cor:CT_single}}\label{app:CT_single_pf}

Analogously to Theorem~\ref{thm:id} we define the polynomial
\begin{align*}
p_{G}(x) &\coloneqq \prod_{i\in\I}\left(x-e^{\lambda_i\tau}\right)^{\tilde{m}_i} \\
	&=  x^{r}+\alpha_{r-1}x^{r-1}+\cdots+\alpha_{1}x+\alpha_{0}.
\end{align*}
Examining \eqref{eq:Jordan_exp_ct} suggests the substitution $\bin{n}{k} = \frac{n^k}{k!}$ in the proof of Theorem~\ref{thm:id}, whose application yields the values $\eta_i \coloneqq e^{\lambda_i \tau}$. Then, the  eigenvalues corresponding to observable eigenmodes may be obtained via $\lambda_i = \log(\eta_i)/\tau$. \qed

\bibliography{MH-bib}

\end{document}